\documentclass[12pt]{article}
\usepackage{amsmath,amsfonts,amssymb,latexsym}
\usepackage[colorlinks, urlcolor=black]{hyperref}
\usepackage{scalerel}
\begin{document}
\baselineskip 15pt

\def \s {\smallskip}
\def \m {\medskip}
\def \b {\bigskip}
\def \n {\noindent}
\def \sn {\s\n}
\def \mn {\m\n}
\def \bn {\b\n}
\def \vn {\vskip-.05in\noindent}
\def \pr {\qquad\qquad}
\def \ol {\overline}
\def \Iff {\Leftrightarrow}
\def \ga {\gamma_n}

\author{Victor Olevskii}
\title{Completion by perturbations}
\date{}
\maketitle

\begin{abstract}
We prove that any non-complete orthonormal system in a Hilbert
space can be transformed into a basis by small perturbations.

\s
\end{abstract}

\b\b

{\bf 1. Introduction.}
Let $H$ be a separable Hilbert space. The standard way to make 
a non-complete system of orthogonal vectors an orthogonal basis 
is adding of the appropriate quantity of missing vectors. 

\s
Our goal here is to show another approach to the problem of 
completion: one can get a basis just by a gentle perturbation 
of the original vectors.

\m\s

{\bf Theorem.}
{\it Let $\{\chi_n\}$, $n\in{\mathbb N}$, be an orthonormal 
system of vectors in $H$. Then there exists an orthonormal 
basis $\{\psi_n\}$ such that
\vskip-.12in

$$\|\psi_n-\chi_n\|\,\xrightarrow[\!\!n\to\infty\!\!\!]{}\,0\,.$$
In addition, all the norms can be made smaller than a given
$\epsilon>0$. 
}

\m\s
Earlier we proved ([O]) such a theorem in a specific case when  
$\{\chi_n\}$ were the indicator functions of the intervals 
$[n, n\!+\!1]$ in $L_2({\mathbb R})$. 

\s
That result answered a question by H. Shapiro whether an 
orthonormal basis can be obtained by translates of a compact 
family of functions. Recall that a system of translates of a 
single function, or a finite family of functions, is never an 
orthonormal basis, nor even a frame, in $L_2({\mathbb R})$ 
([OZ], [CDH]).

\m
It turns out that basically the same method as in [O] allows 
us to perform the result in its general form above.

\newpage

{\bf 2. Proof of Theorem.}
We introduce the following $(n+1)\times (n+1)$ matrices:
$$
A_n=\left(
\begin{array}{ccccc}
1-\ga & -\ga & \ldots & -\ga & \frac{1}{\sqrt{2n}}\\
-\ga & 1-\ga & \ldots & -\ga & \frac{1}{\sqrt{2n}}\\
\vdots & \vdots & \ddots & \vdots & \vdots \\
-\ga & -\ga & \ldots & 1-\ga & \frac{1}{\sqrt{2n}}\\
-\frac{1}{\sqrt{2n}} & -\frac{1}{\sqrt{2n}}
& \ldots & -\frac{1}{\sqrt{2n}} & \frac{1}{\sqrt{2}}\\
\end{array}
\right),
$$
where $\ga=\frac{1}{(2+\sqrt{2})\,n}$.
It is easy to check that the matrices are orthogonal. 

\m
Their key feature is that the upper left $n\times n$ 
submatrix  is close to the identity matrix while the 
lower right element is essentially between 0 and 1.

\m
Our construction below is inspired by Bourgain's paper [B].
We will use the following

\m
{\bf Lemma.} {\it Let $\Psi=\{\psi_n\}$ be a set in a Hilbert 
space. For some $\alpha <1$, suppose there is a set $\Gamma$, 
dense in the unit sphere, such that every $g\in \Gamma$ can 
be approximated, with an error less than $\alpha$, by a linear 
combination of vectors $\psi_n$. Then $\Psi$ is complete.
}

\m
{\it Proof.} If $\Psi$ were not complete, then there would 
be a vector $f$, $\|f\|=1$, orthogonal to span$(\Psi)$. 
Take $g\in \Gamma$ with $\|f-g\| < 1-\alpha$, and find 
$\psi\in$ span$(\Psi)$ with $\|g-\psi\| < \alpha$. 
Then $\|f-\psi\| < 1$ \,contradicts the choice of $f$ 
$\bullet$

\b
Fix a sequence of integers $0<n_1<n_2<...$\,, and split
the original system $\{\chi_n\}$ into blocks consisting of 
the corresponding number of vectors:

\m
\centerline{
$\chi_1^{(1)}=\chi_1, \,\chi_2^{(1)}=\chi_2,
\,...\,, \,\chi_{n_1}^{(1)}=\chi_{n_1}\,; \qquad$
}

\m
\centerline{
$\,\qquad \chi_1^{(2)}=\chi_{n_1+1}, \,\chi_2^{(2)}=\chi_{n_1+2},
\,...\,, \,\chi_{n_2}^{(2)}=\chi_{n_1+n_2}\,;$
}

\m\s\n
etc. Denote the linear span of the $k$-th block by $L^{(k)}$.
The subspaces $L^{(k)}$ are pairwise orthogonal. Define also 
$L_k=L^{(k)}\oplus L^{(k+1)}\oplus ...$\,. 
The orthogonal complement of $L_1$, call it $L_0$, 
is non-empty, otherwise $\{\chi_n\}$ would be a basis.

\m
Let now $\Gamma=\{g_k\}$, $k\in \mathbb N$, be a dense set 
in the unit sphere of $H$, with an additional property:
$g_k \perp L_k$ $\forall k$. Such a set exists since  
$\,\mathbin{\scaleobj{0.9}{\bigcap}}_k\, L_k=\{\overline 0\}$,
i.e. $\mathbin{\scaleobj{0.9}{\bigcup}}_k\, L_k^{\perp}\,$
is dense in $H$.

\m
The desired orthonormal basis will be built up inductively.


\textbf{\textit{Step 1.}} \,Take $n=n_1$, and apply the matrix 
$A_n$ to the orthonormal set of vectors \,$\chi_1^{(1)}$,
$\chi_2^{(1)}$, ...\,, $\chi_n^{(1)}$, and $g^{(1)}=g_1$:
\vskip-.06in

$$
A_n\left(
\begin{array}{c}
\chi_1^{(1)}\\
\vdots \\
\chi_n^{(1)}\\
g^{(1)}\\
\end{array}
\right)
=\left(
\begin{array}{c}
\psi_1^{(1)}\\
\vdots \\
\psi_n^{(1)}\\
h^{(1)}\\
\end{array}
\right).
\quad $$
The obtained vectors \,$\psi_1^{(1)}$, $\psi_2^{(1)}$, ...\,, 
$\psi_n^{(1)}$, and $h^{(1)}$ are also orthonormal. We have:

\s
\centerline{
$\|\psi_j^{(1)} - \chi_j^{(1)}\|^2 =
\|(1-\ga)\chi_j^{(1)} - \ga\sum_{i\ne j}\chi_i^{(1)}
+ \frac{1}{\sqrt{2n}} g^{(1)} - \chi_j^{(1)}\|^2 =$
}

\m\s
\centerline{
$=\|\frac{1}{\sqrt{2n}}g^{(1)}-\ga\sum_{i=1}^n\chi_i^{(1)}\|^2
=\frac{1}{2n} +\ga^2\,n=\frac{1}{2n} (1+\frac{1}{3+2\sqrt{2}})
\,< \frac{1}{n}\,.$
}

\b\n
Since $A_n^{-1} = A_n^{T}$ and, consequently,
$g^{(1)}= \frac{1}{\sqrt{2n}} \sum_{i=1}^n\psi_i^{(1)} +
\frac{\,1}{\sqrt{2}} h^{(1)}$, we also get

\m\s
\centerline{
$\|\sum_{i=1}^n\frac{1}{\sqrt{2n}}\psi_i^{(1)} - g_1\|^2
\,= \,\frac{1}{2}\,.$
}

\b

\textbf{\textit{Step 2.}} Take now $n=n_2$. By definition, 
the vector $g_2\in L_2^{\perp}=L_0\oplus L^{(1)}$. 
The dimension of this subspace is greater than $n_1$, 
so we may assume that $g_2$ does not belong to the linear 
span of $\{\psi_1^{(1)}, ...\,, \psi_{n_1}^{(1)}\}$: 

\m\s\n
\centerline{
$g_2 = \sum_{i=1}^{n_1}
\langle g_2,\psi_i^{(1)}\rangle\,\psi_i^{(1)} + 
\lambda\, g^{(2)},$
}

\m\s\n
where $g^{(2)}$ is orthogonal to all $\{\psi_i^{(1)}\}$,  
$\|g^{(2)}\|=1$, and $0<\lambda\leq 1$.
Apply the operator $A_n$ to the orthonormal set of vectors 
\,$\chi_1^{(2)}$, ...\,, $\chi_n^{(2)}$, and $g^{(2)}$:
\vskip-.06in

$$
A_n\left(
\begin{array}{c}
\chi_1^{(2)}\\
\vdots \\
\chi_n^{(2)}\\
g^{(2)}\\
\end{array}
\right)
=\left(
\begin{array}{c}
\psi_1^{(2)}\\
\vdots \\
\psi_n^{(2)}\\
h^{(2)}\\
\end{array}
\right).
\quad $$
The obtained vectors \,$\psi_1^{(2)}$, ...\,, $\psi_n^{(2)}$, 
and $h^{(2)}$ are again orthonormal. Besides, 
$\psi_j^{(2)}\bot \,\psi_i^{(1)}$ \,$\forall\,i,j$.
Arguing exactly as in Step 1, we have

\m
\centerline{
$\|\psi_j^{(2)} - \chi_j^{(2)}\|^2 < \frac{1}{n}$
}

\n
and

\centerline{
$\|\sum_{i=1}^n\frac{1}{\sqrt{2n}}\psi_i^{(2)}
- g^{(2)}\|^2 = \frac{1}{2}\,.$
}


\n
Due to the decomposition of $g_2$, this equality implies,
for some numbers $\lambda_i^{(k)}$, that

\vskip-.02in
\centerline{
$\|\sum_{k=1}^2\sum_{i=1}^{n_k}\lambda_i^{(k)}\psi_i^{(k)}
- g_2\|^2 \le \frac{1}{2}\,.$
}

\b\s

\textbf{\textit{Step 3 and so forth}} \,are similar to Step 2, 
just $g^{(k)}$ should be taken orthogonal to all 
$\psi_j^{(1)}$, ...\,, $\psi_j^{(k-1)}$, not only to 
$\psi_j^{(k-1)}$.
Finally, we obtain an orthonormal system of functions
$\{\psi_j^{(k)}\}$, \,$j=1, ..., n_k$, $k\in {\mathbb N}$,
such that

\m\s
\centerline{
$\|\psi_j^{(k)}-\chi_j^{(k)}\| < 
\frac{1}{\sqrt{\scaleobj{1.3}{n_k}}}\,,$
}

\m\n
and their linear combinations approximate every $g_k$ up to 
$\frac{1}{\sqrt{2}}$. By the lemma, the system is complete. 
Denote it, in the sequential numbering, by $\Psi=\{\psi_n\}$. 
We have

\vskip-.07in
\centerline{
$\|\psi_n-\chi_n\|\,\xrightarrow[\!\!n\to\infty\!\!\!]{}\,0\,.$
}

\m\n
In addition, taking $n_1 > \frac{1}{\epsilon^2}$, we get the 
estimate $\,\|\psi_n-\chi_n\|<\epsilon\ \,\forall n$ \
$\bullet$

\b

{\bf 3. Remark.} The speed of convergence is optimal, 
in the following sense. Choosing the numbers $\{n_k\}$ growing 
exponentially, one can reach 

\m\n
\centerline{
$\|\psi_n-\chi_n\| \,=\, O(\frac{1}{\sqrt{\scaleobj{1.3}{n}}})$
}

\vskip.09in\n
for $n\to\infty$. 
On the other hand, the norms must satisfy the equality

\m\n
\centerline{
$\sum_n \|\psi_n - \chi_n\|^2 = \infty$\,,}

\m\n
otherwise, according to a well-known theorem of N. Bari,
the orthonormal systems $\{\psi_n\}$ and $\{\chi_n\}$ 
are complete or non-complete simultaneously.

\m
\renewcommand{\refname}{\normalsize\bf \ \ \quad References}

{\footnotesize
\n
{\sc
Dept. of Mathematics, Afeka Tel Aviv, Israel.} \,\ 
{\it E-mail:} {\tt victoro@afeka.ac.il}
}

\end{document}